\documentclass[10pt,letterpaper]{article}
\usepackage{opex3}
\usepackage{color}

\usepackage[utf8]{inputenc}

\usepackage[english]{babel}
\usepackage{amsmath,amssymb,mathrsfs,amsfonts,listings,dsfont,amsthm,units}
\usepackage{xcolor}
\usepackage{cite}
\usepackage{graphicx}

\newcommand{\mR}{\mathbb{R}}
\newcommand{\mC}{\mathbb{C}}

\newcommand{\cD}{\mathcal{D}}

\newcommand{\cF}{\mathcal{F}}

\newcommand{\cO}{\mathcal{O}}
\newcommand{\cP}{\mathcal{P}\,}

\newcommand{\cR}{\mathcal{R}}

\newcommand{\bh}{\ensuremath{\boldsymbol{h}}}

\newcommand{\bxi}{\ensuremath{\boldsymbol{\xi}}}

\newcommand{\D}{\text{d}}
\newcommand{\I}{\text{i}}
\newcommand{\norm}[1]{\left \| #1 \right \|} 

\newcommand{\conj}{\overline}

\newcommand{\NF}{N_{\text{F}}}

\DeclareMathOperator*{\argmin}{argmin}
\DeclareMathOperator*{\sign}{sign}

\begin{document}


\title{Regularized Newton methods for x-ray phase contrast and general imaging problems}

\author{Simon Maretzke,$^{1,*}$ Matthias Bartels,$^2$ Martin Krenkel,$^2$ Tim Salditt$^2$ and Thorsten Hohage$^{1}$}

\address{$^1$Insitute for Numerical and Applied Mathematics, Lotzestr. 16-18, 37083 Goettingen, Germany \\
$^2$Insitute for X-Ray Physics, Friedrich-Hund-Platz 1, 37077 Goettingen, Germany}

\email{$^*$s.maretzke@math.uni-goettingen.de} 



\begin{abstract}
Like many other advanced imaging methods, x-ray phase contrast imaging and tomography require mathematical inversion of the observed data to obtain real-space information. While an accurate forward model describing the generally nonlinear image formation from a given object {to the observations} is often available, explicit inversion formulas are typically not known. Moreover, the measured data might be insufficient for stable image reconstruction, in which case it has to be complemented by suitable \emph{a priori} information. In this work, regularized Newton methods are presented as a general framework for the solution of such ill-posed nonlinear imaging problems. For a proof of principle, the approach is applied to x-ray phase contrast imaging in the near-field propagation regime. Simultaneous recovery of the phase- and amplitude from a single near-field diffraction pattern {without homogeneity constraints} is demonstrated for the first time. The presented methods further permit \emph{all-at-once} phase contrast tomography, i.e.\ {simultaneous} phase retrieval and tomographic inversion. We demonstrate the potential of this approach 
by three-dimensional imaging of a colloidal crystal at $\unit[95]{nm}$ isotropic resolution.
\end{abstract}

\ocis{(100.3190) Inverse problems; (100.3010) Image reconstruction techniques; (100.5070) Phase retrieval; (100.6950) Tomographic image processing; (340.7440) X-ray imaging.} 


\bibliographystyle{osajnl}
\renewcommand\refname{References and links}


\section{Introduction} \label{S1}
     
Lensless coherent diffractive x-ray imaging (CDI) has opened up a new field of high resolution structure analysis beyond the ensemble averaging of conventional x-ray diffraction \cite{Nugent2010coherent,Quiney2010CDIReview,Miao2015XrayCDIReview}.
Typically, lensless x-ray imaging setups are closer to diffraction experiments than to a classical microscope setup with lenses, except that they require
a sufficiently coherent probing wavefront or beam as well as sufficient sampling of the diffraction pattern. In other words, the imaging systems is essentially based on free
space propagation between object and detector. Depending on whether the data is recorded in the optical near-field or far-field, the propagation is modeled by the Fresnel propagator 
or a Fourier transform (Fraunhofer far-field regime), respectively. As in conventional diffraction, each detector pixel carries information about all object pixels. Therefore, if the data is modeled in the detector plane, high spatial resolution can only be achieved in terms of the average structure. Contrarily, if the data is {\it inverted} to reconstruct the object, the individual real space configuration is depicted, beyond an ensemble average (or more precisely an average over the entire illuminated volume). Importantly,  the image 
formation is still very similar to a plain diffraction experiment, while the data analysis is not. Data Modeling for example by least-square fitting is replaced by image reconstruction. Model formulation is replaced by the formulation of \emph{a priori} knowledge (constraints), which are required to compensate for the missing information on the phase of the diffraction field, and hence to achieve a unique solution.  Accordingly, the difference is brought about by the inversion of the diffraction process: by solving the inverse 
problem of (non-crystallographic) diffraction, we obtain access to the individual configuration instead of the average sizes and correlation lengths in the object. 
It is for this reason that iterative algorithms \cite{MiaoNature1999CDIFirstExp,Bauschke2002PhaseRetr,Elser2003ItProjPhaseRetr,PaganinXRay,Marchesini2007EvalPhaseRetrAlg}, being to date the engine of CDI, have received so much attention. 
Iteratively cycling between the detector and object planes, they feed in both measured data and additional \emph{a priori}  information on the solution.

CDI uses \textit{a priori} knowledge, for example related to the object support or its optical constants (positivity, pure phase contrast). 
Diversity in the data generated by illuminating  the same object pixels by different wavefronts may be exploited by ptychographic algorithms  \cite{Faulkner2004Ptycho,Rodenburg2007Ptycho,Thibault2008Ptycho,Maiden2009ImprovedPtycho,Thibault2013ImagingStateMixt}. 
In general, however, one is interested in finding experimental settings allowing for robust phase retrieval with least invasive constraints, dose and 
accumulation times. In many applications for example, one is interested in reconstructing from single recordings, without scanning or multiple exposures. 
At the same time, the object may be composed of several materials, both with phase and amplitude contrast.  
It is therefore often of advantage to perform the imaging in the optical near-field rather than far-field.
Here, phase information is much more directly encoded in the  diffraction pattern, by ways of interference between 
diffracted and reference waves. Recently, improved  uniqueness results have been presented 
for near-field imaging using the theory of entire functions \cite{Maretzke2015IP}. In fact, contrary to a common 
belief, see e.g. \cite{Jonas2004TwoMeasUniquePhaseRetr,Burvall2011TwoPlanes}, 
measurements  at only one distance are sufficient to determine both the real and the imaginary part of any 
compactly supported refractive index distribution. These results are much stronger than those for the Fraunhofer regime, as obtained e.g.\ by \cite{Klibanov20062dUniquePurePhase}.
Loosly speaking, the near-field case can thus be expected to require less information to complement the measurements  than the far-field case.
In practice, however, many  iterative algorithms which have been adapted to the near-field case \cite{Giewekemeyer2011CellSketch,Stockmar2013NearFieldPtycho,Robisch2015NearFieldPtycho} require similar data diversity or constraints.

 Notwithstanding the merits of many highly performing iterative algorithms, it is therefore necessary to broaden the perspective of the phase retrieval field.
Here, we present iteratively regularized Gauss-Newton methods (IRGNM) \cite{Bakushinskii1992IRNM} as an alternative approach to phase retrieval and other imaging problems. In this method, each iterate is computed to provide an optimal compromise between agreement with the measured data and additional constraints on the basis of a local linearization of contrast formation, as we discuss further below. {The approach is related to the regularized gradient descent methods for phase reconstruction proposed in \cite{Davidoiu2011FrechetPhaseRetr}, but promises improved convergence owing to the Newton-like solution of linearized subproblems.}
{We apply the general IRGNM framework to near-field phase contrast with x-rays,} 
for the reasons mentioned above, and since recent progress in propagation imaging has narrowed the gap in resolution compared to CDI.
In fact, using highly divergent and coherent quasi-spherical wavefronts, x-ray imaging in the optical near-field has been recently demonstrated down to 20 nm resolution \cite{Bartels2015}.

 The scope of the present work is two-fold. Firstly, we give a concise overview of iteratively regularized Gauss-Newton methods in view of x-ray imaging, summarizing and explaining recent mathematical literature in this field for an applied audience.  Secondly, we demonstrate the performance of this 
approach in solving the phase problem on the level of real state-of-the-art experimental data.
In particular, we apply an IRGNM approach to three-dimensional (3d) imaging, 
i.e.\@ we show how the method can be used to perform phase retrieval and tomographic reconstruction simultaneously. {This strategy has been argued to enable improved accuracy compared to  splitting the reconstruction into phase retrieval problems for each angle to recover the fields in the object 
plane and a subsequent inversion of the Radon transform \cite{Chapman2006,Barty2008ceramicfoam,Ruhlandt2014,MyMaster}.}
We show that IRGNM approaches offer significant flexibility in treating various experimental setups and different \textit{a priori} information.
In the long run, we also expect advantages owing to the fact that the noise characteristics of the data can be suitably accounted for in this framework.
For the most important example of Poisson data, Newton-type regularization methods with Kullback-Leibler-type data fidelity terms have already been proposed \cite{Hohage2013}.

 The manuscript is organized as follows: $\S$\ref{S2} introduces IRGNM in view of image reconstructions problems from the principal idea to practical implementation. In $\S$\ref{S3}, the framework is applied to (2d) near-field phase contrast with hard x-rays, imaging the phase shifts and absorption induced by a nano-structured test pattern. A Kaczmarz-type IRGNM suitable for tomographic  imaging is presented in  $\S$\ref{S4} and applied to resolve the structure of a colloidal micro-crystal. 

\section{Regularized Newton methods for imaging} \label{S2}

\subsection{Basic approach} \label{SS2.1}

We consider an abstract imaging system of the form
\begin{equation}
 I^{\text{obs}} = F ( f^\dagger )  + \boldsymbol{\varepsilon}. \label{eq:AbsImageProb}
\end{equation}
Here $I^{\text{obs}} \in Y $ denotes some observable intensity data, given by the image of the unknown object $f^\dagger\in X$ (e.g. a spatially varying refractive index or scatterer positions) under a known forward operator $F$ and superimposed measurement errors $\boldsymbol{\varepsilon}$. Note that $\boldsymbol{\varepsilon}$ may depend on $F ( f^\dagger ) $ as is the case for Poisson data. For many models, an explicit inversion formula for $F$ is not known. Moreover, even if the inverse $F^{-1}$ is available, it is often not continuous and would thus amplify errors $\boldsymbol{\varepsilon}$ by large magnitudes if applied  directly to the data $I^{\text{obs}}$. Such \emph{ill-posedness} of the reconstruction despite uniqueness is well-established e.g.\ for computed tomography \cite{Natterer}. We seek a method to stably recover $f^\dagger$ from an ill-posed problem of the form \eqref{eq:AbsImageProb}.

Notably, the operator $F: X \to Y$ modeling the imaging system is nonlinear in general, for example whenever a phase retrieval problem described by a squared modulus operation is involved. Nevertheless, reasonable results can often be achieved using a linearization of contrast formation, as given for instance by the contrast transfer function in electron microscopy or coherent x-ray imaging \cite{Erickson1971_CTF_TEM,Wade1992_CTF_TEM,Guigay1977CTF,Cloetens1999}.

Mathematically, such first order approximations are justified by the Fr\'echet differentiability of the forward operator $F$, i.e.\ for any $f$ there exists a bounded linear map $F'[f]$ such that $\lim_{h\rightarrow 0}(F ( f + h) - F(f) - F'[f]h) / \norm{h} = 0$. A natural approach to solve \eqref{eq:AbsImageProb} is then given by Newton-type iterations 
\begin{equation}
 f_{k+1} = f_k  + F'[f_k]^{-1} \left(  I^{\text{obs}} - F ( f_k ) \right). \label{eq:NaiveNewton}
\end{equation}
As opposed to methods based on a \emph{static} linearization of contrast formation, the linearizations in \eqref{eq:NaiveNewton} are computed about the current iterate $f_k$ and thereby account for moderate nonlinearity.

{However, iterations of the form \eqref{eq:NaiveNewton} are often neither feasible nor desirable for imaging because the linearized problems - just like the nonlinear equation \eqref{eq:AbsImageProb} - are typically ill-posed. Hence, the solution of \eqref{eq:NaiveNewton} is unstable and 
may not even exist.} Again, one may think of ambiguities in phase retrieval problems.
A remedy is given by iteratively regularized Gauss-Newton methods (IRGNM) as first proposed by \cite{Bakushinskii1992IRNM}, corresponding to Tikhonov regularization of the Newton steps:
\begin{align}
 f_{k+1} = \argmin_{f \in X}   \left\| F(f_{k}) + F'[f_k](f - f_k) - I^{\text{obs}} \right\|_{Y}^2&    \nonumber \\ + \; \alpha_k  \left\| f - f_{0} \right\|_X^2&  \label{eq:IRGNM} 
\end{align}
Here, $\| \cdot \|_{Y}$ and $\| \cdot \|_{X}$ denote the norms in Hilbert spaces $X$ and $Y$, $f_0\in X$ is the initial guess and $\alpha_k > 0$ is a regularization parameter.
In this setting, it can be shown that \eqref{eq:IRGNM} always has a unique solution given by \cite{Hanke1996Regularization}
\begin{align}
 f_{k+1} = f_k + \left(F'[f_k]^\ast F'[f_k] + \alpha_k\right)^{-1} \big( \;\;\;\; &F'[f_k]^\ast ( I^{\text{obs}} - F ( f_k ) ) \nonumber \\ 
									    + \; &\alpha_k ( f_0 - f_k ) \; \big). \label{eq:NewtonStep}
\end{align}
$F'[f_k]^\ast$ denotes the adjoint of the linear map $F'[f_k]: X \to~Y$. In \eqref{eq:NewtonStep}, only the inverse of the selfadjoint positive-definite operator $F'[f_k]^\ast F'[f_k] + \alpha_k$ has to be computed, which is \emph{bounded} according to the estimate
\begin{equation}
 \norm{ \left( F'[f_k]^\ast F'[f_k] + \alpha_k\right)^{-1}} \leq \frac 1 {\alpha_k}. \label{eq:RegStability}
\end{equation}
Accordingly, the iterate $f_{k+1}$ depends continuously on $I^{\text{obs}}$, i.e.\ the impact of data errors on the reconstruction is \emph{regularized}.

As the IRGNM is  based on linearizations of the imaging operator, yet iteratively updated, the approach is best suited for weakly or moderately nonlinear problems. Formally, convergence of the method to $f^\dagger$ for $\boldsymbol{\varepsilon} \to 0$ can indeed be shown given bounds on the nonlinearity of $F$ along with suitably chosen $\alpha_k$ and $f_0$ \cite{Blaschke1997IRGNMConv,Kaltenbacher2008ItRegNonlin}.

\subsection{Parameter choice and constraints} \label{SS2.2}

 The first term on the right hand side of \eqref{eq:IRGNM} measures the agreement of the object with the observed data $ I^{\text{obs}}$ based on the current linearization. To achieve competitive results, the choice of the norm $\norm{\cdot}_Y$ should reflect the statistical properties of the data errors $\boldsymbol{\varepsilon}$, e.g. by taking the negative log-likelihood of the measured signal $I^{\text{obs}}$. For additive Gaussian white noise, this consideration leads to the choice of the standard $L^2$-norm, i.e.\
  \begin{equation}
  \norm{I - I^{\text{obs}}}_Y := \int  \left| I - I^{\text{obs}} \right|^2 \; \D x.
 \end{equation}
 For Poisson noise, the resulting data fidelity term is the Kullback-Leibler-divergence. This distance measure can be implemented in the framework of generalized Newton methods as demonstrated by \cite{Hohage2013}. Within the IRGNM, a quadratic approximation about its minimum may be used as a norm
 \begin{equation}
  \norm{I - I^{\text{obs}}}_Y := \norm{\frac{ I - I^{\text{obs}} }{  \max( I_0, I^{\text{obs}} )^{\frac 1 2 } }  }_{L^2},  
 \end{equation}
where $I_0 > 0$ is a regularizing parameter.

On the right hand side of \eqref{eq:IRGNM}, the data residual is balanced with the regularization term $\alpha_k  \left\| f - f_{0} \right\|_X^2$ bounding the deviation from the initial guess $f_0$. The regularization parameter weights the different contributions: if $\alpha_k$ is very small, there is essentially no regularization and the norm bound \eqref{eq:RegStability} diverges, allowing for large amplifications of the data error. If $\alpha_k$ is chosen too large on the other hand, the Newton iterate computed via \eqref{eq:NewtonStep} need not have much to do with the actual image to be reconstructed. A good choice of $\alpha_k$ must thus balance data- and approximation errors. One possible strategy is the following:
\vspace{.5em}
\begin{itemize}
 \item Choose $\alpha_0$ to approximately balance the norms in \eqref{eq:IRGNM}, for example by setting $\alpha_0 \sim \norm{F'[f_0] F'[f_0]^\ast I^{\text{obs}}}_Y^2 / \norm{F'[f_0]^\ast I^{\text{obs}}}_X^2$. \vspace{.25em}
  \item Reduce $\alpha_k$ by a constant factor, e.g. $\alpha_{k+1} / \alpha_k = \frac 2 3$. \vspace{.25em}
  \item Stop at the first  $k $ s.t.\ $s_k := \norm{F(f_{k}) -  I^{\text{obs}}}_Y \leq \tau \norm{\boldsymbol{\varepsilon}}_Y$, i.e.\ when the residual attains $\tau \geq 1$ times the noise level. \vspace{.5em}
\end{itemize}
The stopping criterion, known as Morozov's discrepancy principle \cite{Morozov1966Discrepancy}, requires good knowledge of the magnitude of data errors. When this is not available, noise level-free stopping rules need to be applied, see e.g.\ \cite{Hanke1996Regularization,Hansen2010DiscreteIP,Kindermann2011NoiseFreeMethods}. Here, the principal idea is to make $s_k$ as small as possible while limiting heuristic measures for the impact of data errors such as the object norm $\norm{f_k - f_0}_X$ or the inverse regularization parameter $\frac 1 {\alpha_k}$.

By the choice of the norm $\norm{\cdot }_X$, we may define desirable properties of the object $f$ to be reconstructed. Choosing the standard $L^2$-norm here prevents isolated spikes, promoting more evenly distributed values. Typical images are expected to be of higher regularity, for example being composed of smoothly varying areas bounded by sharp edges. This may be exploited to obtain higher robustness to high-frequency errors by regularizing with Sobolev norms
\begin{equation}
 \norm{f}_X := \norm{ (1 + \bxi^2)^{\frac s 2}  \cF(f)( \bxi ) }_{L^2} \label{eq:SobolevReg}
\end{equation}
where $\cF$ denotes the Fourier transform and $\bxi$ the frequency coordinates. The exponent $s \geq 0$ tunes the degree of smoothness.

Beyond smoothness it is often desirable to impose additional constraints, which corresponds to restricting the set of admissible solutions $f_{k} \in C \subset X$. Prominent examples are real-valuedness, support constraints or positivity. {Geometrically, the former two types are represented by linear subspaces $C \subset X$. Imposing these constraints in the IRGNM simply amounts to substituting $F'[f_k]^\ast$ with $ \cP F'[f_k]^\ast $ in \eqref{eq:NewtonStep}, where $\cP: X \to C$ is the orthogonal projection onto~$C$.} Positivity, on the other hand, is a nonlinear convex constraint. It may be included within a generalized Newton framework via a nonsmooth regularization term. In practice, this amounts to solving the minimization problem \eqref{eq:IRGNM} restricted to $f\in C$, as can be done using semismooth Newton methods \cite{Hintermuller2002SSNewton}. One approach to approximate sign constraints within the IRGNM framework of the present work lies in supplementing \eqref{eq:IRGNM} with the penalty term
\begin{equation}
 \gamma \norm{ \min( 0, f_k ) - \min( 0, \sign(f_k) ) ( f - f_k) }_{L^2}^2, \label{eq:PosTerm}
\end{equation}
which tends to correct negative values of $f_k$ in the subsequent iterate. The coefficient $\gamma > 0$ determines the weight of the constraint and should be comparable to $\alpha_0$. To achieve strict positivity, one may let $\gamma \to \infty$ at constant $\alpha_k$ in the final iterates.

For numerical implementation of the IRGNM, all that needs to be done is to exchange the imaging operator $F$ and its derivative $F'$ in \eqref{eq:IRGNM} by suitable discrete approximations. The norms $\norm{\cdot}_{X_{\text{dis}}}$  and $\norm{\cdot}_{Y_{\text{dis}}}$ in the discretized object- and image spaces $X_{\text{dis}}= \mR^{N_X}$, $Y_{\text{dis}} = \mR^{N_Y} $ are characterized by their Gramians $\textbf{G}_X, \textbf{G}_Y$ with respect to the  Euclidean norm, i.e.\
\begin{equation}
 \norm{\boldsymbol f }_{X_{\text{dis}}}^2 = \boldsymbol f^{\text{T}} \textbf{G}_X  \boldsymbol f, \;\;  \norm{\boldsymbol I }_{Y_{\text{dis}}}^2 = \boldsymbol I^{\text{T}} \textbf{G}_Y  \boldsymbol I \label{eq:Gramians}
\end{equation}
where $\boldsymbol f^{\text{T}}$ denotes the transpose. The adjoint of the Fr\'echet derivative $F'_{\text{dis}}[\boldsymbol f]: X_{\text{dis}} \to Y_{\text{dis}}$ can be implemented via
\begin{equation}
 F'_{\text{dis}}[\boldsymbol f]^\ast = \textbf{G}_X^{-1} F'_{\text{dis}}[\boldsymbol f]^{\text{T}} \textbf{G}_Y. \label{eq:DiscAdj}
\end{equation}
The  Hermitean positive-definite linear problem in the Newton step \eqref{eq:IRGNM} can be solved efficiently by a conjugate gradients (CG) method. A major advantage of this IRGNM-CG scheme is that only the \emph{forward} maps $F_{\text{dis}}$, $F'_{\text{dis}}[\boldsymbol f] $ and $F'_{\text{dis}}[\boldsymbol f]^{\text{T}}$ need to be evaluated, which are often easy to implement. Moreover, note that the matrices $F'_{\text{dis}}[\boldsymbol f]$ usually need not (and should not) be set up explicitly to compute the matrix-vector products $\bh \mapsto F'_{\text{dis}}[\boldsymbol f]^{(\text{T})} \bh$.

\subsection{Newton-Kaczmarz methods} \label{SS2.3}

The Newton-type method presented so far is an \emph{all-at-once} approach to image reconstruction: linearizations of the problem \eqref{eq:AbsImageProb} are solved in each iteration, incorporating all constraints and the complete measured data $I^{\text{obs}}$. In view of 3d or even 4d tomographic or time-resolved imaging data, this results in computations with huge arrays, posing numerical challenges. {For \emph{linear} problems, a well-known remedy is given by Kazcmarz-type methods \cite{Kaczmarz1937ART} such as the Algebraic Reconstruction technique in computed tomography \cite{Gordon1970ART, Natterer}, which cyclically solve small under-determined subproblems.} Indeed, many nonlinear imaging problems also lend themselves to a natural separation into different subproblems so that \eqref{eq:AbsImageProb} becomes
\begin{equation}
 (I^{\text{obs}}_1, \ldots, I^{\text{obs}}_p) = (F_1, \ldots, F_p) (f) + \boldsymbol{\varepsilon}.
\end{equation}
For problems of this form, \cite{Burger2006NewtonKaczmarz} proposed regularized Newton-Kaczmarz methods, where in each iteration only one of the linearized subproblems
$ F_{j_k}(f_k) + F'_{j_k}[f_k](f_{k+1} - f_k) = I^{\text{obs}}_{j_k}$ is considered. A Kaczmarz-type equivalent of the IRGNM-update \eqref{eq:IRGNM} takes the form
\begin{align}
 f_{k+1} = \argmin_{f \in X}   &\left\| F_{j_k}(f_{k}) + F_{j_k}'[f_k](f - f_k) - I^{\text{obs}}_{j_k} \right\|_{Y_{j_k}}^2   \nonumber \\
 + \; \alpha_k \big( \beta_k  &\left\| f - f_{0} \right\|_X^2 + \; (1-\beta_k)  \left\| f - f_{k} \right\|_X^2 \big)  \label{eq:Kaczmarz-Newton} 
\end{align}
The additional regularization term bounds the deviations from the preceding iterate $f_{k}$ according to the weights $\beta_k\in [0;1]$, ensuring that the reconstruction does not change too much within one Newton iteration. The sequence $(j_k)$ determines the processing order of the different subproblems which may be adjusted to the requirements of the particular imaging problem.

\section{Application to propagation-based phase contrast} \label{S3}

\begin{figure}[htbp]
\centering
\includegraphics[width = 0.9\textwidth]{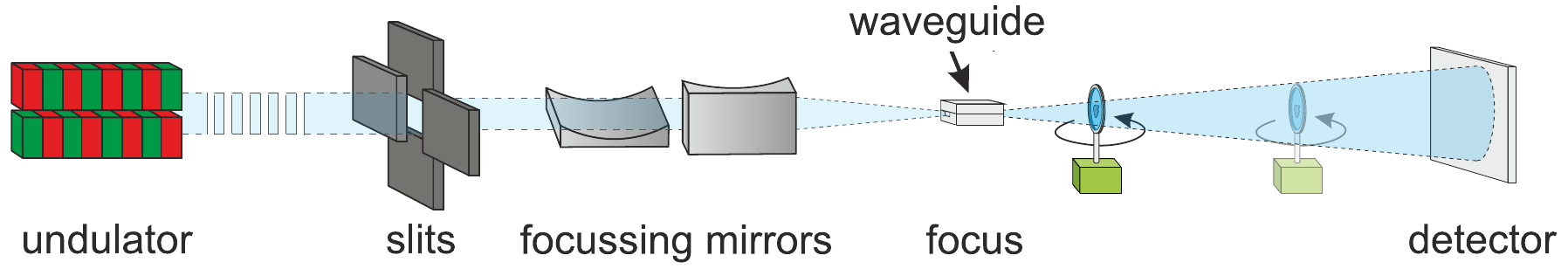}
\caption{Schematic setup for propagation-based phase contrast imaging with hard x-rays: GINIX at P10 beamline, DESY \cite{Krenkel2014BCAandCTF,Salditt2015GINIX}. Quasi-monocromatic x-rays are focused onto a waveguide, illuminating a downstream object by a cone beam emanating from this coherent quasi-point source. Rotation of the sample allows for tomographic measurements. \label{fig:1}}
\end{figure}

Next, we apply the regularized Newton framework presented in section \ref{S2} to propagation-based phase contrast x-ray imaging. An exemplary experimental setup is sketched in Fig. \ref{fig:1}: {quasi-monochromatic undulator radiation is focused by a pair of elliptical mirrors onto an x-ray waveguide. The exit of the waveguide serves as a quasi point source illuminating an object, which is placed at a distance of several mm behind the waveguide.} The resulting diffraction patterns (holograms) are recorded by a detector at about 5m distance behind the focal plane,
capturing the entire cone beam emanating from the  waveguide with the sample induced interference pattern.  
The imaging data presented in this work has been recorded at the GINIX endstation at P10 beamline, DESY, Hamburg, described in  \cite{Kalbfleisch2011GINIX,Salditt2015GINIX}. For details on the waveguide system and comparable high resolution propagation imaging results, see also \cite{Bartels2015}.

For the imaging model, a fully coherent illumination of the sample by a plane wave is assumed. Note that we consider an effective parallel beam geometry equivalent to the experimental cone-beam setup. As in previous studies of cone beam x-ray propagation imaging \cite{Pogany1997noninterferometric,Mayo2002quantitative}, we make use
of the simple variable transformation, accounting for the geometric magnification, mapping the spherical beam illumination to an effective plane wave case.
For the present imaging system, the data treatment is detailed by \cite{Bartels2015}. Moreover, we assume the object to be sufficiently thin and weak for the wave field to remain spherical and for the transmission to be well described by geometrical optics (projection approximation). This is typically well-satisfied for hard x-ray imaging \cite{Mayo2002quantitative,PaganinXRay}. Under these assumptions, the intensity at the detector can be modeled as
\begin{equation}
 I  = F (f) := \left|\cD \left(\exp( -\I f )\right) \right|^2 \label{eq:PCIModel},
\end{equation}
where $f = \phi - \I \mu/2$ parametrizes the phase shifts $\phi$ and attenuation $\mu$ induced by the specimen, i.e.\ its image. $\cD$ denotes the Fresnel propagator, implementing free-space propagation of the transmitted parabolic wave field from sample onto the detector:
\begin{equation}
 \cD ( \psi) := \cF^{-1} \left( \exp\left( - \frac{\I\pi \bxi ^2}{\NF}   \right) \cF( \psi ) (\bxi )  \right)
\end{equation}
$\NF = a^2 k / (2\pi d)$ is the Fresnel number of the setup, where $a$ is the size of a feature corresponding to unit length in the frequency vector $\bxi$ (typically the size of a pixel), $d$ is the defocus distance and $k$ is the x-ray wavenumber.

From \eqref{eq:PCIModel}, we see that phase contrast imaging is nonlinear in general. Nevertheless, the model can be linearized for weak phase shifts and attenuation $\phi, \mu$, yielding the contrast transfer function (CTF) \cite{Guigay1977CTF,Turner2004FormulaWeakAbsSlowlyVarPhase}
\begin{align}
 \cF( F(f) - 1) \approx  2 &\sin \left( \frac{ \pi \bxi ^2}{\NF} \right) \cF(\phi) (\bxi)  \nonumber \\
 - &\cos \left( \frac{ \pi \bxi ^2}{\NF} \right) \cF(\mu) (\bxi)  \label{eq:CTF}
\end{align}
Assuming vanishing attenuation $\mu \approx 0$ or single-material objects $\mu \propto \phi$, this formula can be directly inverted for image reconstruction. This approach yields satisfactory reconstructions especially if data from multiple defocus distances is available \cite{Cloetens1999,Krenkel2014BCAandCTF}.

From the CTF solution \eqref{eq:CTF}, it can be seen that regularized Newton methods are well-suited for phase contrast imaging in two respects: firstly, the zeros of the oscillatory prefactors make the image reconstruction problem ill-posed even in the simplest case of a weak non-absorbing object so that regularization is required. Secondly, the relative success of linear CTF-based methods indicate that the nonlinearity of the imaging method is sufficiently weak in a large regime of experimental setups.

For suitable spaces $X,Y$, the imaging operator $F: X\to Y$ defined by \eqref{eq:PCIModel} can be shown to be Fr\'echet differentiable. {Using linearity of $\cD$ and the expansions $\exp(x+h) = \exp(x) ( 1 + h) + \cO(h^2)$ and $| x + h|^2 = |x|^2 + 2 \Re(\conj x h) + \cO(h^2)$, we obtain the first-order Taylor approximation
\begin{equation}
 F(f+h) = F(f) +  2 \Re \left(-\I \overline{\cD \left(\exp( -\I f )\right)} \cD \left( \exp( -\I f ) h \right) \right) + \cO( \norm{h}_X^2). \label{eq:PCIDerivativeDerivation}
\end{equation}
This allows to read off the Fr\'echet derivative given by the linear term in $h$:}
\begin{equation}
 F'[f]h = 2 \Im \left( \overline{\cD \left(\exp( -\I f )\right)} \cD \left( \exp( -\I f ) h \right) \right). \label{eq:PCIDerivative}
\end{equation}
Here the overbar denotes complex conjugation and $\Im$ the pointwise imaginary part. The special case $f= 0$ in \eqref{eq:PCIDerivative} reproduces the CTF. Implementing \eqref{eq:PCIModel} and \eqref{eq:PCIDerivative} for discrete images in the spaces $X := \mC^{M_X \times N_X}$ and  $Y := \mR^{M_Y \times N_Y}$, using suitably padded fast Fourier transforms for the discrete Fresnel propagator, the
reconstruction problem of propagation-based phase contrast imaging can be solved by the IRGNM of section \ref{S2}.

\begin{figure}[htbp]
\centering
\includegraphics[width = .9\textwidth]{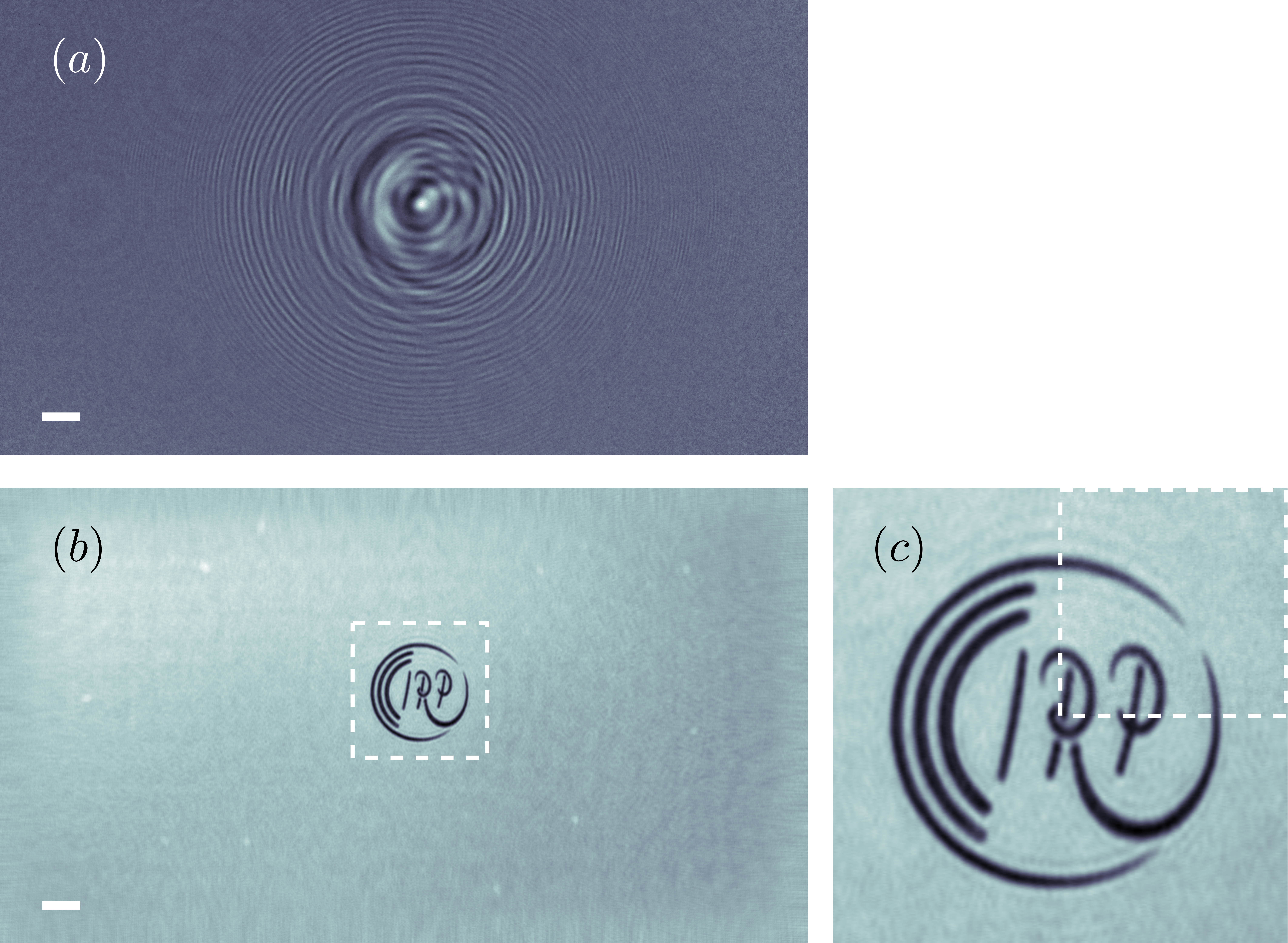}
\caption{X-ray phase contrast imaging of a nano-structured object using the IRGNM reconstruction algorithm. The test pattern (institute logo) was defined in a thin gold film by focused ion beam milling. $(a)$ Diffraction pattern recorded at GINIX
endstation, P10 beamline, DESY \cite{Salditt2015GINIX} (flat-field corrected, maximum photon flux per pixel $\approx 3400$). $(b)$~Reconstructed phase image of the entire field of view, assuming a fixed ratio of $0.21$ between absorption $\mu$ and phase shifts $\phi$. $(c)$ Magnification of $(b)$ in the framed region around the logo. For comparison, the dashed inset shows the corresponding part of the phase map reconstructed by direct inversion of the CTF \eqref{eq:CTF} via the methods of \cite{Krenkel2014BCAandCTF}. Scale bars: $\unit[2]{\mu m}$ in the effective geometry. Fringe artifacts on the upside of the logo are due to diffraction fringes leaving the field of view in $(a)$. The phase shifts $\phi \approx 0.20$ induced by the object correspond to a thickness of the gold layer of about $\unit[100]{nm}$.\label{fig:2}}
\end{figure}

As a proof of concept, we consider the diffraction pattern of an IRP logo, engraved into a thin gold film. The data has been recorded on the aforementioned GINIX setup at an energy of $E= \unit[7.9]{keV}$ is shown in Fig. \ref{fig:2}$(a)$. The plotted intensity data has been normalized by the corresponding flat-field, i.e. the image of the empty beam, in order to meet the assumption of plane wave illumination as justified by the analysis of \cite{Homann2015Validity}. The Fresnel number is $\NF = 1.77 \cdot 10^{-4}$ at an effective pixel size of $a = \unit[21.7]{nm}$ of the $1080 \times 1920$-sized images and a maximum flux per pixel of $\approx 3400$ photons.

For the Newton reconstruction, we choose the Sobolev norm in \eqref{eq:SobolevReg} with an exponent $s = 0.5$ as a regularization term and an $L^2$-norm for the data residual. The regularization parameters $\alpha_k$ are chosen according to the procedure outlined in section \ref{SS2.2}. The algorithm is stopped as soon as the reduction of the residual $s_k- s_{k-1}$ achieved in the $k$-th iteration falls below $1\,\%$ of the maximum decrease $\max_{j<k} s_j- s_{j-1}$.  For the initial guess, we simply take $f_0 = 0$. We exploit that the specimen is made of gold by prescribing its characteristic between absorption $\mu$ and phase shifts $\phi$ of $c_{\mu/\phi} \approx 0.21$ \cite{Henke1993TypicalDeltaBeta}. The constraint is imposed by substituting $f = (1- \frac{\I } 2 c_{\mu/\phi} ) \phi$ and reconstructing the real-valued $\phi$. Apart from this, neither positivity nor a support constraint is assumed. The latter means that the phase shifts $\phi$ are computed within the entire field of view of $1080 \times 1920$ pixels without any oversampling.

Convergence of the IRGNM is reached after 12 Newton steps, corresponding to a total of 242 CG-iterations. The resulting phase map is plotted in Fig. \ref{fig:2}$(b)$. Despite the lack of oversampling, it can be seen that both the object and the background come out quite clean except for some stripes near the boundary due to the applied smooth replicate padding. Occasional white spots and the smooth background variations, on the other hand, may be  attributed to some dirt in the imaging optics and physical variations of the thickness of the gold layer, respectively. The magnification of the reconstruction in a region around the IRP-logo depicted in Fig. \ref{fig:2}$(c)$ confirms the high uniformity of the recovered phase shifts in the regions within and without the logo, reproducing the binary test pattern. Slight fringe artifacts can be identified in the vertical direction, especially on the upper edge of the logo.
{These may be explained by incompleteness of the data. Indeed, it can be seen from Fig. \ref{fig:2}$(a)$ that parts of the diffraction pattern lie outside of the recorded field of view. The reconstruction of the missing object information is based on \emph{a priori} constraints only, which naturally gives rise to artifacts.}

It should be noted that the imaging setting considered so far still lies well within the regime of applicability of reconstruction method based on the CTF \eqref{eq:CTF}: for once, the phase shifts within the logo of $\approx \unit[0.20]{rad}$, induced by the $\unit[100]{nm}$-thick gold layer of refractivity $k \delta \approx \unit[1.96]{\mu m ^{-1}}$  \cite{Henke1993TypicalDeltaBeta}, are sufficiently weak for a global linearization of contrast formation to be reasonably accurate. Secondly, the assumption of a single material, i.e.\ of proportional fields $\phi$ and $\mu$, allows for a direct inversion of \eqref{eq:CTF}. In the considered example, the latter approach indeed yields results of comparable quality as shown by the dashed inset in Fig. \ref{fig:2}$(c)$. {If $\phi$ and $\mu$ are considered as independent variables and only a single diffraction pattern is available, directly inverting \eqref{eq:CTF} is impossible. This suggests that the imaging problem is non-unique in this case.} On the contrary, theoretical investigation has recently shown that a unique recovery is indeed feasible if $\mu$ and $\phi$ are compactly supported \cite{Maretzke2015IP}.
The mathematical reason is that compactness of the support translates into correlations in Fourier space, which can be used to disentangle $\phi$ and $\mu$ in \eqref{eq:CTF}. By its ability to incorporate \emph{a priori} constraints, the IRGNM approach may thus enable reconstructions beyond the limitations of direct inversion methods.
\begin{figure}[htbp]
\centering
\includegraphics[width = .8\textwidth]{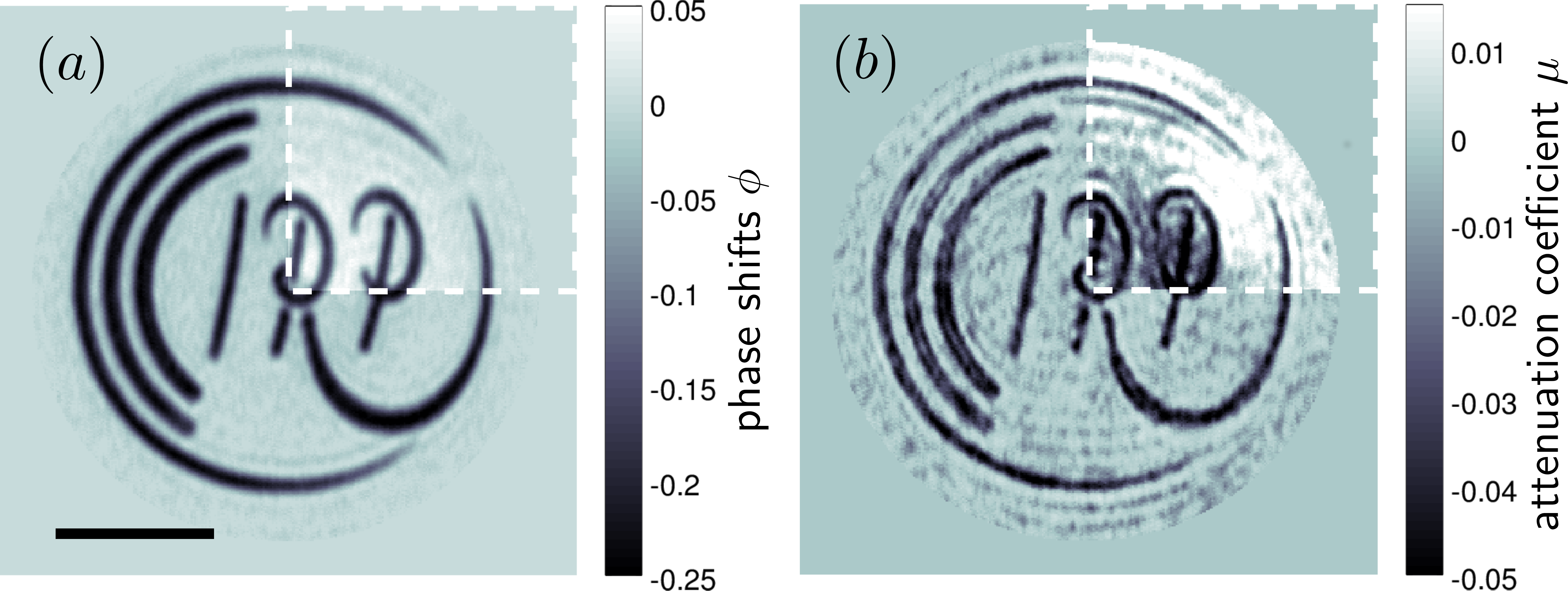}
\caption{Simultaneous IRGNM reconstruction of phase shifts $\phi$ and absorption $\mu$ from the data in Fig. \ref{fig:2}$(a)$ without assuming a fixed ratio $\mu / \phi$. Negative values of $\mu$ and $\phi$ indicate missing material in the gold film. The circular support visible in the phase- and absorption maps and negativity of $\mu$ and $ \phi$ has been imposed as a constraint. All other parameters are retained as in the computation of Fig. \ref{fig:2}$(b)$.
{The dashed inset in the absorption image shows the reconstruction without 
negativity constraint in the IRGNM and demonstrates that simultaneous recovery 
tends to introduce low-frequency artifacts.}
These are effectively suppressed by exploiting physical \emph{a priori} knowledge on the sign of phase shifts and absorption. Scale bar: $\unit[2]{\mu m}$.\label{fig:3}}
\end{figure}

For a first experimental verification of a simultaneous recovery of both phase and absorption from a single hologram, we repeat the IRGNM reconstruction using the same parameters but without imposing a coupling of $\phi$ and $\mu$. Following the uniqueness analysis, we impose a loose circular support constraint around the logo.
{Note that this induces a strong oversampling in the data by a factor of $\approx 17$, owing to the considerably reduced number of object pixels on which $\phi$ and $\mu$ have to be reconstructed.}
Moreover, negativity of $\mu$ and $\phi$ (the test pattern represents \emph{missing} material!) is imposed in the IRGNM via the penalty term approach outlined in section \ref{SS2.2}.
The recovered phase and absorption obtained after 10 Newton steps (199 CG-iterations) is plotted in Fig. \ref{fig:3}. For comparison, dashed insets show the corresponding image parts obtained from a second reconstruction without negativity constraints.


Both the phase shifts $\phi$ in \ref{fig:3}$(a)$ and the attenuation coefficient $\mu$ plotted in \ref{fig:3}$(b)$ clearly represent the object shape, where the magnitudes roughly reproduce the material-specific ratio of $\mu / \phi \approx 0.21$. The only visible artifacts are those attributed to the missing fringes, which have been observed previously. In the IRGNM-reconstruction without negativity constraints, additional low-frequency errors appear especially in the recovered attenuation, as can be seen from the central dark spot surrounded by a bright halo in the inset of Fig. \ref{fig:3}$(b)$. Simulations as well as preliminary analytical studies indicate that the susceptibility to these halo artifacts is characteristic of simultaneous retrieval of phase and absorption. 
As observed in the present reconstruction, however, they seem to be suppressed very effectively by imposing suitable sign constraints.
Owing to the flexibility of exploiting physical \emph{a priori} knowledge in the IRGNM, we thus indeed achieve an almost perfect simultaneous recovery of $\mu$ and $\phi$ from a single hologram - which has to date been considered impossible \cite{Nugent2007TwoPlanesPhaseVortex,Burvall2011TwoPlanes}.
\begin{figure}[htbp]
\centering
\includegraphics[width = .8\textwidth]{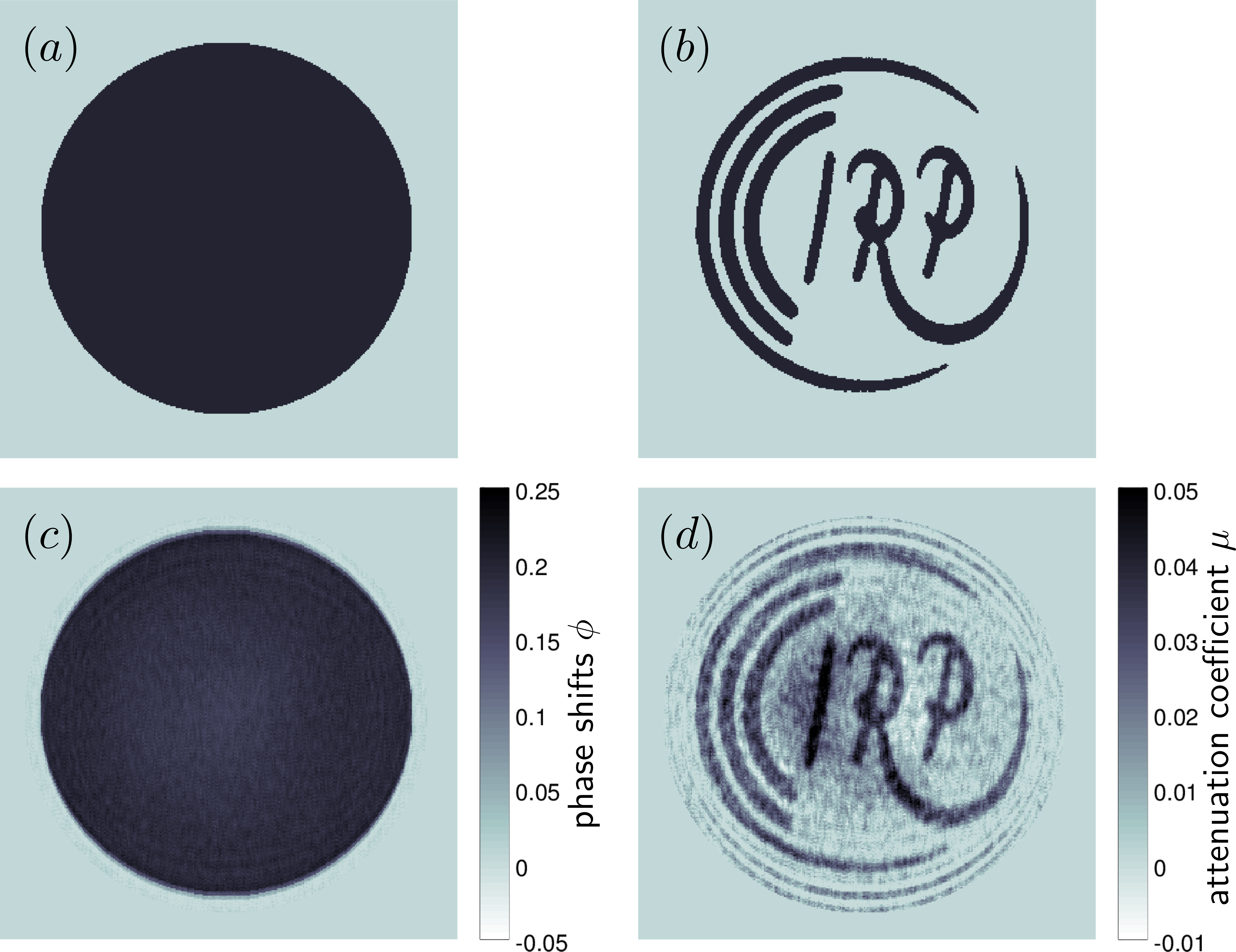}
{\caption{Supplementary IRGNM reconstruction of a non-homogeneous object from a simulated near-field hologram with $2\,\%$ Gaussian white noise. Constraints, reconstruction- and setup parameters are chosen as for the experimental data in Fig. \ref{fig:3}. $(a)$ and $(b)$: Exact object composed of an absorbing logo structure with $\phi = 0.2$, $\mu = 0.04$ (gold) embedded into a non-absorbing disc that induces the same phase shifts $\phi = 0.2$. $(c)$ and $(d)$: Recovered phase and absorption images. Except for a slight halo, features in $\phi$ and $\mu$ are correctly identified, revealing the hidden logo structure by the enabled absorption contrast.  \label{fig:3.2}}}
\end{figure}

{Although proportionality of $\mu$ and $\phi$ is not imposed in the reconstruction in Fig.~\ref{fig:3}, the surprising quality might still be owing to some implicit preference of the IRGNM for homogeneous objects. We test this by repeating the reconstruction for a simulated object that is composed of two different materials. To this end, a binarized version of the recovered IRP logo with constant $\phi = 0.2$ and $\mu = 0.04$ 
is embedded into a purely phase-shifting disc with $\phi = 0.2$ and $\mu = 0$.
All setup- and reconstruction parameters are chosen exactly as in Fig.~\ref{fig:3}. A single synthetic hologram (not shown) is simulated by mapping the exact object with the forward operator and superimposing Gaussian white noise with standard deviation $\sigma = 0.02$.
The IRGNM reconstruction of  $\phi$ and $\mu$ from this data is shown in Fig.~\ref{fig:3.2} along with the exact solutions.

Good agreement with the exact object is found both in phase and absorption. The different object features are correctly attributed to the $\phi$- and $\mu$-components in Figs.~\ref{fig:3.2}$(c)$ and $(d)$  and no artifacts except for the known slight halo and fringe structures are visible. In particular, the absorption contrast enabled by the simultaneous recovery clearly reveals the logo structure despite the low signal-to-noise-ratio in $\mu$. The logo is not represented in the phase image and thus likely invisible to any reconstruction method that assumes proportionality of $\mu$ and $\phi$.}

\section{Phase contrast tomography of a colloidal crystal} \label{S4}

By rotating the specimen within the incident beam in the imaging setup of Fig. \ref{fig:1}, propagation-based phase contrast can be extended to a tomographic imaging method, capable of resolving three-dimensional variations of the complex refractive index $n = 1 - \delta + \I \beta$. Within the geometrical optics approximation of section \ref{S3}, the imprinted phase- and absorption images $\phi$ and $\mu$ are proportional to the projection of $\delta $ and $\beta$, respectively, along the optical axis. For a tomographic incident angle $\theta$, we have
\begin{equation}
  \phi_\theta - \I \mu_\theta/2 = k \int  (\delta_\theta - \I \beta_\theta) \; \D z = k \cR ( \delta - \I \beta)(\theta)   \label{eq:ProjRadon}
\end{equation}
where $\delta_\theta$ and $\beta_\theta$ denote the fields $\delta$ and $\beta$ in a rotated coordinate system and $\cR$ is the two-dimensional Radon transform  mapping onto the rotated line integrals. Combining \eqref{eq:PCIModel} and \eqref{eq:ProjRadon}, we obtain a forward operator $F_{\text{tomo}}$ mapping the object function $f_{\text{tomo}} = \delta - \I \beta$ onto the corresponding intensity data $I_{\text{tomo}}$ under \emph{all} of the different incident angles:
\begin{equation}
 I_{\text{tomo}}  = F_{\text{tomo}} (  f_{\text{tomo}} ) := \left|\cD \left(\exp( -\I k \cR ( f_{\text{tomo}} )   )\right) \right|^2 \label{eq:PCTFwOp}
\end{equation}
Since $\cR$ is linear and bounded, $F_{\text{tomo}} $ is Fr\'echet differentiable with a derivative of similar form as \eqref{eq:PCIDerivative}.

Implementing \eqref{eq:PCTFwOp} in a regularized Newton-type framework corresponds to an \emph{all-at-once} approach to phase contrast tomography, as phase- and tomographic reconstruction are not carried out in subsequent steps but simultaneously. In particular, this implies that tomographic correlations between images under different incident angles are incorporated already in the phase retrieval step. Similarly as in coherent diffractive imaging \cite{Chapman2006,Barty2008ceramicfoam}, this has been shown to improve stability and accuracy of the reconstruction \cite{Ruhlandt2014,MyMaster}. Unfortunately, a numerical inversion of $F_{\text{tomo}}$ by an IRGNM is not viable for high resolution data sets due to memory constraints and high computational costs of the evaluation of the Radon transform.

As a remedy, we propose to solve the reconstruction problem of phase contrast tomography by a regularized Newton-Kaczmarz methods as presented in section \ref{SS2.3}. The tomographic setting suggests a decomposition of the forward operator $F_{\text{tomo}} = (F_1, \ldots, F_d)$, corresponding to small sub-data sets of diffraction patterns from only very few incident angles. Each Newton step \eqref{eq:Kaczmarz-Newton} then only requires evaluations of comparably small Radon transforms and processing of feasible chunks of data. At the same time, tomographic consistency is exploited by directly reconstructing a three-dimensional object instead of single projections. Furthermore, the $F_j$ may be chosen such that the most strongly correlated images for only slightly differing incident angles are processed simultaneously. 

We implement the Newton-Kaczmarz step \eqref{eq:Kaczmarz-Newton} for the imaging operator \eqref{eq:PCTFwOp}, using the simple parameter choice $\alpha_k = \alpha_0$ and $\beta_k = \beta_0 \ll 1$, where $\alpha_0$ is estimated according to section \ref{SS2.2}. Firstly, this regularization imposes the current reconstruction as a strong prior for the next iteration. Secondly, it also yields a favorable condition number for the linear Newton steps so that CG-methods typically converge after 3-5 iterations.

\begin{figure}[htbp]
\centering
\includegraphics[width = .8\textwidth]{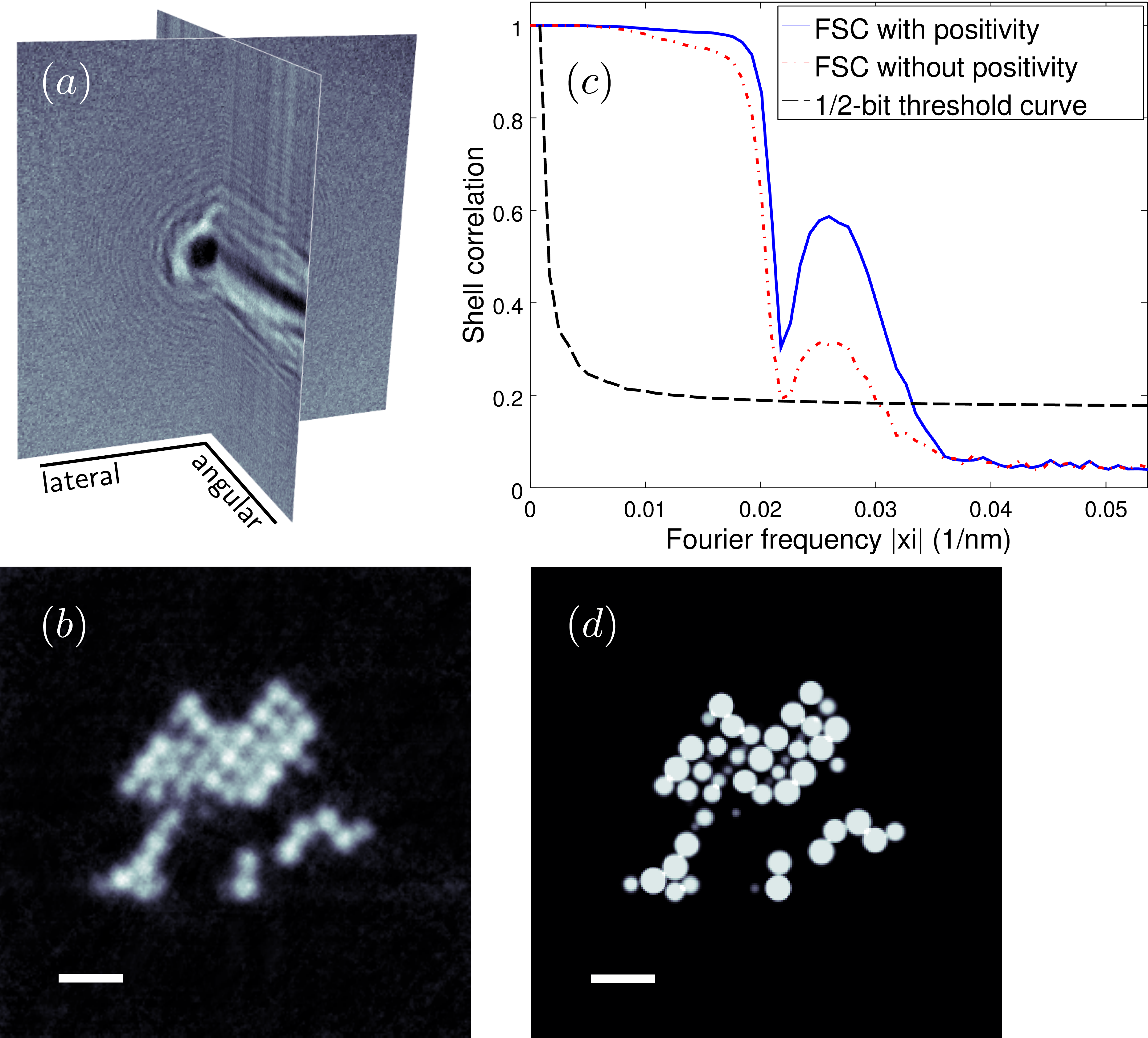}
\caption{X-ray phase contrast tomography of a colloidal micro-crystal composed of $\unit[415]{nm}$-size polystyrene beads, reconstructed by a regularized Newton-Kaczmarz method. The maximum photon flux is $770$ per pixel corresponding to an absorbed dose $\leq \unit[110]{kGy}$. $(a)$ Slice visualization of the tomographic data given by $1024\times 1024$ flat-field corrected holograms measured under 249 incident angles between $0$ and $\unit[173]{^\circ}$ at GINIX endstation, P10 beamline, DESY \cite{Salditt2015GINIX}. 
$(b)$ Central slice of the reconstructed $256^3$ voxel volume, plotting the increment $\delta$ of the refractive index $n = 1-\delta + \I  \beta$ which is proportional to electron density. \emph{A priori} constraints were positivity $\delta \geq 0 $ and vanishing absorption $\beta = 0$. The scale bar is $\unit[1]{\mu m }$. $(c)$ Fourier shell correlation (FSC) computed for reconstructions from complementary data sets of 125 and 124 incident angles (green curve). For comparison, the blue curve shows the result without positivity constraint and the red dashed one plots the $1/2$-bit threshold curve \cite{Heel2005OneHalfBitThres}. The intersection between the green- and the red curve indicates an achieved resolution of $\unit[95]{nm}$. $(d)$~Binary representation of the slice in $(b)$, determined by deconvolving the reconstructed $\delta$ with the form factor of a $\unit[415]{nm}$-sized homogeneous sphere after Gaussian filtering.\label{fig:4}}
\end{figure}

For a proof of concept, we consider phase contrast images of a colloidal crystal composed of $\unit[415]{nm}$-sized polystyrene beads, imaged at an energy $E = \unit[7.9]{keV}$ in the experimental setup of Fig. \ref{fig:1}. The tomographic data set is composed of 249 diffraction patterns of size $1024 \times 1024$ acquired at a maximum photon flux of $770$ per pixel under incident angles $\theta$ between $0$ and $\unit[173]{^\circ}$. The raw data correction included an alignment of the tomographic projections to the common center of mass and iterative reprojections, following the approach of \cite{Xu_NatMat_2015}.
The photon flux corresponds to an upper bound of $\unit[110]{kGy}$ for the absorbed dose over the total tomographic data acquisition. The flat-field corrected intensities are visualized in Fig. \ref{fig:4}$(a)$. In the effective geometry, the Fresnel number is $\NF  = 2.41 \cdot 10^{-4}$ at an effective pixel size of $\unit[29.3]{nm}$. 
  
From the fringes in the data set in Fig. \ref{fig:4}$(a)$ as well as preliminary reconstructions, one can infer the object location in the center of the field of view. This is exploited by restricting the reconstruction to a central $256^3$ voxel cube, imposing a loose 3d-support constraint. Moreover, the hydrocarbon composition of the polystyrene spheres and 
the photon energy allow using $\beta = 0$ as a constraint (non-absorbing object), i.e.\ to reconstruct only the refractive increment $\delta$. Unfortunately, a Sobolev-type regularization term \eqref{eq:SobolevReg} would bottleneck the proposed Newton-Kaczmarz-scheme as the required FFTs in the Gramian would constitute the only $\cO(N^3\log N )$ operations in the Newton updates. We therefore recur to simple $L^2$-regularization and $L^2$-data fidelity terms, corresponding to identity Gramians $\textbf{G}_X, \textbf{G}_Y$ in \eqref{eq:Gramians}. {However, positivity of $\delta$ is imposed by supplementing the Kaczmarz-Newton step \eqref{eq:Kaczmarz-Newton} with the penalty term in \eqref{eq:PosTerm}. The weights of the penalty terms are set to $\gamma = \alpha_0$ and $\beta_k = 0.001$.} For the initial guess, we simply choose $f_0 = 0$. The forward operator $F_{\text{tomo}}$ is decomposed such that angular ``wedges'' of six adjacent diffraction patterns are processed simultaneously in each Newton step. The number of iterations is determined such that each hologram is visited twice where the processing order of the wedges is random to minimize directional bias in the reconstructed object.

The reconstruction runs on a simple laptop with  Intel i7 CPU and 8 Gigabytes of RAM, computing for about 20 minutes. Each of the Newton steps takes four CG-iterations, which means that the total solution - processing every incident angle twice - requires only 10 forward and adjoint operations of the full imaging operator $F_{\text{tomo}}$. Figure \ref{fig:4}$(b)$ shows the central slice of the reconstructed $\delta$ perpendicular to the tomographic axis. It can be seen that the binary character of the specimen and the spherical shapes of the beads are well-resolved, being clearly distinguishable both from one another and from the background. Notably, no artifacts caused by angular undersampling or due to the missing wedge of $6^\circ$ can be identified.

We estimate the resolution by computing the Fourier shell correlation (FSC) \cite{Heel1987FSC} of two auxillary reconstructions from complementary sets of 125 and 124 incident angles, respectively. In Fig. \ref{fig:4}$(c)$, the resulting correlation is compared to the $1/2$-bit threshold curve as proposed by \cite{Heel2005OneHalfBitThres}. For comparison, the FSC is also computed for a reconstruction without positivity constraint. The beneficial effect of the latter can be seen from the gap between the blue and green curves in Fig. \ref{fig:4}$(c)$, showing a significantly higher correlation especially in the resolution-critical part around the intersections with the threshold curve. The intersections indicate a resolution of about $\unit[95]{nm}$ and $ \unit[105]{nm}$ (half of the corresponding Fourier wavelengths) for the reconstructions with and without positivity constraint, respectively. Notably, the local minima of the FSC curves at $|\xi| \approx \unit[0.022]{nm^{-1}}$ are no artifacts but neatly coincide with the first order zero of the form factor of $\unit[415]{nm}$-sized spheres. This emphasizes the sensitivity of the imaging method to structural features.

According to the slice plotted in Fig. \ref{fig:4}$(b)$, the positions of the polystyrene spheres and thus the crystal structure could be determined in principle by visual inspection of the reconstructed $\delta$. In order to determine the colloid locations numerically and independent of an operator, we smoothen the recovered object by a Gaussian of $\unit[95]{nm}$ FWHM to reduce the impact of noise before deconvolving with the form factor of a homogeneous $\unit[415]{nm}$-sized sphere (with regularization around the zeros in Fourier space). This procedure results in Gaussian peaks centered at the positions of polystyrene beads. By computing the maxima of this field using quadratic interpolation about the peaks, we thus obtain a representation of the crystal structure.

\begin{figure}[htbp]
\centering
\includegraphics[width = .6\textwidth]{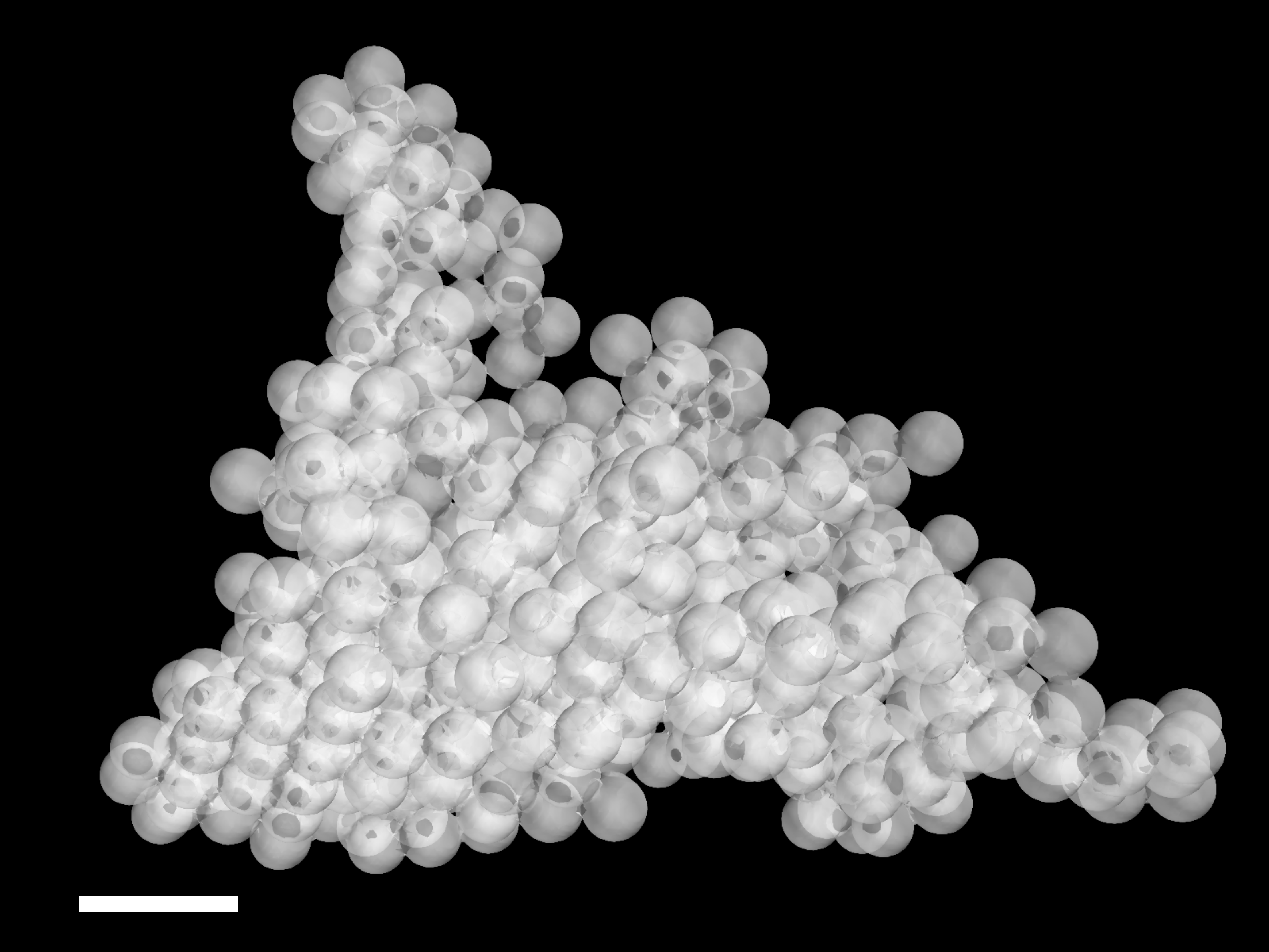}
\caption{3d-rendering of a colloidal crystal of polystyrene nano-beads reconstructed by Newton-based phase contrast tomography. The underlying binary object, a slice of which is shown in Fig. \ref{fig:4}$(d)$, has been determined by deconvolving the smoothened reconstruction of the refraction parameter $\delta$ with the form factor of a homogeneous sphere of diameter $\unit[415]{nm}$. The coordinates of the colloid sites are provided in Data File 1. Scale bar: $\unit[1]{\mu m}$. \label{fig:5}}
\end{figure}

By the simple deconvolution procedure outlined above, 448 colloid positions are identified.
Convolving the obtained Dirac delta-array of bead locations with the ideal $\unit[415]{nm}$-sized sphere yields a binary representation of the colloidal crystal. A slice of this binarization corresponding to Fig. \ref{fig:4}$(b)$ is plotted in  \ref{fig:4}$(d)$. Figure~\ref{fig:5} shows a 3d-rendering of the determined colloidal crystal. The corresponding site coordinates of the colloids, i.e. the centers of the spherical beads, are provided in Data File 1. The imaged micro-crystal contains regions of (approximately) hexagonal close-packing as well as cubic and amorphous regions, induced by an interplay of bulk- and surface effects.

\section{Conclusions}

In this work, we have presented iteratively regularized Gauss-Newton methods (IRGNM) as a generic approach to solve nonlinear ill-posed image reconstruction problems. The principal idea is to reconstruct an unknown object by iteratively inverting linearizations of a known imaging model, which describes contrast formation from object to observable data. In order to compensate for missing information, the iterates are computed to provide an optimal compromise between the measurements and additional \emph{a priori} information on the unknown object. The IRGNM approach differs from well-known alternating-projection-type algorithms typically used in CDI in that it exploits differentiability and \emph{simultaneously} processes constraints and observed data. This promises improved convergence.

By applying regularized Newton methods to near-field phase constrast imaging, both in 2d and in a tomographic 3d setting, we have demonstrated their flexibility in treating different experimental setups. The reconstruction of a nano-structured test pattern shows that IRGNM constitute a reasonable generalization of direct inversion methods based on a \emph{global} linearization of contrast formation. Owing to its potential to incorporate moderate nonlinearity and \emph{a priori} knowledge, e.g.\ on support or positivity, the approach permits faithful reconstructions beyond the scope of such direct methods. This has been demonstrated by the simultaneous recovery of magnitude and phase of the test object from a single diffraction pattern {without assuming proportionality of phase shifts and absorption.
The validity of the reconstruction is supported by results for a simulated object composed of two materials, which are correctly identified by the algorithm.
However, numerical simulations also indicate that the stability of the simultaneous recovery deteriorates significantly with increasing Fresnel number. Thus, it is possibly only feasible for deeply holographic near field data as considered in this work.}

A further benefit of the presented regularized Newton framework is that it requires only a \emph{forward} model for the imaging setup and no explicit knowledge of (approximate) inverses. This enables IRGNM reconstructions also for complicated multi-staged imaging setups, possibly including the influence of various experimental parameters. In the considered example of phase contrast tomography this flexibility permits to directly recover a probed 3d-object from the complete tomogram of diffraction patterns instead of subsequently performing (2d) phase retrieval and tomographic backprojection. We have demonstrated the potential of this approach by imaging a colloidal crystal composed of $\unit[415]{nm}$-sized polystyrene beads at $\unit[95]{nm}$ resolution, reconstructed by an efficient regularized Newton-Kaczmarz method.

X-ray imaging techniques currently reach out for novel frontiers, ranging from  new contrast modalities 
in material science and dose reduction in biomedical analysis to in-operando studies and investigations of highest spatial and temporal resolution by FEL pulses. 
To enable these goals, the mathematical modeling underlying image reconstruction has to become increasingly accurate in order not to bottleneck the achievable contrast and resolution: apart from nonlinearity of image formation and Poissonian noise statistics, future imaging models may need to better account for partial coherence, 
non-uniform illumination,  mechanical vibrations and movements, misalignment or geometrical aberrations, just to name a few. Owing to their  flexibility, we are convinced that regularized Newton methods, as presented in this work, may greatly foster such developments and thus contribute to a bring lensless x-ray techniques as well as other cutting-edge imaging methods to the level of quantitative structural measurements.

\section*{Acknowledgments}

We thank Janne-Mieke Meijer and Andrei Petukhov for providing us with a colloidal crystal film some time ago, from which we
could later prepare the micro-crystal used here as a test object. MB acknowledges John Miao for helpful discussions regarding the alignment of the colloidal data set. Support by the German Research Foundation DFG through project C02 of the Collaborative Research Center SFB 755 - Nanoscale
Photonic Imaging is gratefully acknowledged.

\end{document}